\documentclass[11pt]{amsart}
\usepackage[margin=1.1in]{geometry}
\usepackage{amsmath,amssymb,amsthm,mathtools}
\usepackage{enumitem}
\usepackage{tikz}
\usepackage[hidelinks]{hyperref}

\newtheorem{theorem}{Theorem}[section]
\newtheorem{lemma}[theorem]{Lemma}
\newtheorem{proposition}[theorem]{Proposition}
\newtheorem{corollary}[theorem]{Corollary}
\newtheorem*{theoremA}{Theorem A}
\newtheorem*{theoremB}{Theorem B}
\newtheorem*{theoremC}{Theorem C}
\newtheorem*{theoremD}{Theorem D}
\theoremstyle{definition}
\newtheorem{definition}[theorem]{Definition}
\newtheorem{example}[theorem]{Example}
\newtheorem{problem}[theorem]{Problem}
\theoremstyle{remark}
\newtheorem{remark}[theorem]{Remark}

\newcommand{\Max}{\operatorname{Max}}
\newcommand{\Jac}{J}
\newcommand{\G}{\mathcal{G}}
\newcommand{\TG}{\mathfrak{T}}
\newcommand{\Pk}{\mathcal{P}_k}
\newcommand{\F}{\mathcal{F}}
\newcommand{\Z}{\mathbb{Z}}
\newcommand{\FF}{\mathbb{F}}
\newcommand{\Bool}[1]{\FF_2^{\,#1}}
\newcommand{\co}[1]{{#1}^{c}}
\newcommand{\rad}{\operatorname{rad}}
\DeclareMathOperator{\Aut}{Aut}

\begin{document}

\title[Domination and maximal intersecting families]{Domination and maximal intersecting families\\
in comaximal graphs of commutative rings}

\author{Elham Mehdi-Nezhad}
\address{Department of Mathematics and Applied Mathematics, University of the Western Cape,
Private Bag X17, Bellville 7535, South Africa}
\email{emehdinezhad@uwc.ac.za}


\begin{abstract}
The comaximal graph $\G(R)$ of a commutative ring $R$ has as vertices the nonunits outside the
Jacobson radical, two being adjacent if they generate $R$ as an ideal. For rings with $k\ge4$
maximal ideals, only the bounds $2\le\gamma(\G(R))\le k$ were known, to our knowledge. We prove
$\gamma=\gamma_t=\gamma_c=k$ whenever $3\le k<\infty$, without finiteness or Artinian hypotheses,
classify the minimum dominating sets, and determine the domatic number. Independent dominating sets
correspond to maximal intersecting families of subsets of $\{1,\dots,k\}$; we obtain a closed
formula for the independent domination number of every finite non-local ring, which can be
exponentially larger than the domination number. The proofs rest on the class sizes of the
blow-up description of $\G(R)$ due to Moconja and Petrovi\'c, which we determine for every ring with
finitely many maximal ideals, finite or infinite; with them, domination problems on $\G(R)$ become
weighted extremal problems on subsets of $\{1,\dots,k\}$, and the formula for the independent
domination number follows from the multiplicative form of these sizes. We also treat rings with
infinitely many maximal ideals, graph polynomials, and lattices whose maximal ideals are prime.
\end{abstract}

\maketitle

\section{Introduction}\label{sec:intro}

Throughout, $R$ is a commutative ring with $1\neq 0$, and $U(R)$, $\Max(R)$ and $\Jac(R)$ denote its
units, its maximal ideals and its Jacobson radical. We write $\Bool k=\FF_2\times\dots\times\FF_2$
($k$ factors) for the product ring, not the field of order $2^k$.

The \emph{comaximal graph} $\Gamma(R)$ of Sharma and Bhatwadekar \cite{SB95} has vertex set $R$,
with $a\neq b$ adjacent if $Ra+Rb=R$. Units are adjacent to all vertices and elements of $\Jac(R)$
only to units \cite[Lemma~2.1]{MSSY08}, so one studies
\[
  \G(R)=\Gamma(R)\bigl[R\setminus(U(R)\cup\Jac(R))\bigr],
\]
denoted $\Gamma_2(R)\setminus\Jac(R)$ in \cite{MSSY08}, $CG_J(R)$ in \cite{MR15} and
$\Gamma_2'(R)$ in \cite{SRD21}. This graph is connected of diameter at most $3$ \cite{MSSY08}; its
girth, radius, center and chromatic number were determined by Moconja and Petrovi\'c \cite{MP11};
see also \cite{Sa14,GKJ24}, and \cite{Ra26c,Ra26d} for the spectrum and the vertex connectivity of
$\G(\Z_n)$.

Put $\Lambda(x)=\{M\in\Max(R):x\in M\}$. If $\Max(R)=\{M_1,\dots,M_k\}$, the \emph{type} of $x$ is
$\tau(x)=\{i:x\in M_i\}$. The vertices of $\G(R)$ are the elements whose type is a nonempty proper
subset of $[k]=\{1,\dots,k\}$, and two vertices are adjacent exactly when their types are disjoint
(Lemma~\ref{lem:adjacency}).

For rings with $k\ge3$ maximal ideals, Mehdi-Nezhad and Rahimi \cite{MR15} found a clique
dominating set of size $k$, so that $2\le\gamma(\G(R))\le k$, and proved $\gamma=3$ for $k=3$; the
comparison with \cite{MR15} is made in detail below. Sinha,
Rao and Davvaz \cite{SRD21} proved $2\le\gamma(\G(\Z_n))\le r$, where $r$ is the number of prime
divisors of $n$, and a similar bound for $\Z_n\times\Z_m$; see also \cite{SR20}. Rather
\cite{Ra26a,Ra26b} computed the independence, independent domination and domination polynomials of
$\Gamma(\Z_n)$ when $n$ has at most three prime divisors. For the comaximal ideal graph, see
\cite{DM16,KK20}. The exact value of $\gamma$ for $k\ge4$, the minimum dominating sets, the domatic
number and the independent domination number were not known. We determine all of them; our main
results are the following.

\begin{theoremA}[Theorems~\ref{thm:gamma-k3}, \ref{thm:gamma-k2}, \ref{thm:infinite}]
Let $R$ be a non-local commutative ring.
\begin{enumerate}[label=\textup{(\alph*)}]
\item If $3\le|\Max(R)|=k<\infty$, then
  $\gamma(\G(R))=\gamma_t(\G(R))=\gamma_c(\G(R))=k$.
\item If $|\Max(R)|=2$, then $\gamma(\G(R))=1$ if $R\cong\FF_2\times F$ for a field $F$, and
  $\gamma(\G(R))=2$ otherwise.
\item If $\Max(R)$ is infinite, then $\G(R)$ has no finite dominating set.
\end{enumerate}
\end{theoremA}

Part (a) needs no finiteness or Artinian hypothesis; the equalities $\omega(\G(R))=\chi(\G(R))=k$
are due to Moconja and Petrovi\'c \cite{MP11}. Part (b) shows that \cite[Corollary~3]{SRD21}, which
states $\gamma=2$ for every ring with two maximal ideals, needs the additional hypothesis
$R\not\cong\FF_2\times F$; for example $\G(\Z_6)=K_{1,2}$. In (c) the cardinal $\gamma(\G(R))$ need not be $|\Max(R)|$
(Example~\ref{ex:boolean-infinite}). Theorem~\ref{thm:gamma-sets} describes all minimum dominating
sets; the cases $k=3$ and $k\ge4$ differ.

\begin{theoremB}[Theorem~\ref{thm:domatic}, Corollary~\ref{cor:Zn}]
Let $R$ have $k\ge3$ maximal ideals, and let $m$ be the minimum over $i$ of the number of elements
lying in $M_i$ and in no other maximal ideal, that is,
\[
  m=\min_i\, |\Jac(R)|\prod_{j\neq i}\bigl|(R/M_j)^{\times}\bigr| .
\]
If $m$ is finite, then $d(\G(R))=m+1$. In particular, if $n$ has at least three prime divisors, then
$d(\G(\Z_n))=1+\varphi(n)/(P-1)$, where $P$ is the largest prime divisor of $n$.
\end{theoremB}

For example $d(\G(\Z_{30}))=3$ and $d(\G(\Z_{105}))=9$, although $k=3$ in both cases.

\begin{theoremC}[Theorems~\ref{thm:indep-char}, \ref{thm:i-closed}, Corollary~\ref{cor:i-rings}]
Let $R$ have $k\ge2$ maximal ideals. The independent dominating sets of $\G(R)$ are exactly the
sets $\{x:\tau(x)\in\F\}$, where $\F$ is a maximal intersecting family of nonempty proper subsets
of $[k]$. If $R$ is finite with residue fields of orders $q_1,\dots,q_k$, then
\[
  i(\G(R)) = \frac{|\Jac(R)|}{2}\sum_{\emptyset\neq A\subsetneq[k]}
  \min\Bigl\{\prod_{j\in A}(q_j-1),\ \prod_{j\notin A}(q_j-1)\Bigr\}.
\]
In particular $i(\G(\Bool k))=2^{k-1}-1$, whereas $\gamma(\G(\Bool k))=k$ for $k\ge3$.
\end{theoremC}

The number of maximal intersecting families on a $k$-set is at least doubly exponential in $k$
(Section~\ref{sec:independent}; see also \cite{BMMV13}). The formula depends on the multiplicative form
of the class sizes \eqref{eq:weights-intro} below and fails for other weights
(Remark~\ref{rem:i-fails}).

\begin{theoremD}[Theorem~\ref{thm:polys}, Corollary~\ref{cor:full-graph-polys}]
For every finite non-local ring $R$, the independence, independent domination and domination
polynomials of $\G(R)$ and of $\Gamma(R)$ are sums over families of nonempty proper subsets of
$[k]$, with summands determined by $|\Jac(R)|$ and $q_1,\dots,q_k$.
\end{theoremD}

These formulas contain those of \cite{Ra26a,Ra26b} for $n$ with at most three prime divisors.

We now describe the method. By \cite[Proposition~2.3, Theorem~2.4]{MP11} every nonempty proper subset of $[k]$ is a type, so
$\G(R)$ is a blow-up of the disjointness graph on these subsets, the Boolean graph \cite{GWY15},
which we call the \emph{type graph} $\TG_k$: the elements of a given type form an independent set,
completely joined to the elements of every disjoint type. We first determine the sizes of these
classes (Lemma~\ref{lem:realization}): the class of type $A$ has cardinality
\begin{equation}\label{eq:weights-intro}
  w_R(A)=|\Jac(R)|\cdot\prod_{i\notin A}\bigl|(R/M_i)^{\times}\bigr|
\end{equation}
for every ring with finitely many maximal ideals, finite or not. The blow-up alone does not
determine the domination parameters: $\gamma$ for $k=2$, the domatic number and $i$ all depend on
these sizes. With them, domination questions on $\G(R)$ become weighted extremal problems on
subsets of $[k]$. Two features of $\TG_k$ are used
throughout. The type $[k]\setminus\{i\}$ is adjacent only to $\{i\}$, so every dominating set meets
each of $k$ disjoint unions of two classes (Lemma~\ref{lem:key}); this controls $\gamma$, the
minimum dominating sets and the domatic number. A set of vertices is independent if and only if
its types pairwise intersect, so independent dominating sets correspond to maximal intersecting
families \cite{EKR61}, and $i(\G(R))$ is the minimum weight of such a family. Such a family contains
exactly one set from each pair $\{A,A^c\}$, so its weight is at least the sum of the smaller weights
over all pairs. For the weights \eqref{eq:weights-intro}, $\log w_R(A)$ is an affine function of the
indicator vector of $[k]\setminus A$; after a small generic perturbation of this function, a lighter
set from each pair can be chosen so that these sets form an intersecting family, and the bound is
attained (Theorem~\ref{thm:i-closed}). The main work of the paper lies in the
classification of the minimum dominating sets, the determination of the domatic number, and the
closed formula for the independent domination number (Theorems~\ref{thm:gamma-sets},
\ref{thm:domatic} and~\ref{thm:i-closed}).

Except for Theorem~A(b),(c), the proofs use only the weighted blow-up and the form of
\eqref{eq:weights-intro}. They therefore apply to comaximal graphs of lattices \cite{AK14}
(Section~\ref{sec:lattices}). For bounded lattices whose maximal ideals are prime, the analogues of
Theorems~A and~B and of the first part of Theorem~C hold. The formula for $i$ fails in general: for
divisor lattices the class sizes do not decrease with the type, and the minimum is attained instead
by the family of all sets containing a fixed index (Proposition~\ref{prop:i-increasing}), by a theorem on product
measures of intersecting families. Without primeness Theorem~A(a) fails: for the diamond $M_3$ one has
$\gamma=1<3=k$.

This paper develops the comaximal-graph part of \cite{MR15}; the ideal-based zero-divisor graph is
not considered here. Theorem~\ref{thm:gamma-k3} turns \cite[Theorem~3.4, Corollary~3.6]{MR15} into
$\gamma=\gamma_t=\gamma_c=k$ for every ring with $3\le k<\infty$, and Theorem~\ref{thm:gamma-sets}
describes all minimum dominating sets. The bounds $2\le\gamma\le|\Max(R)|$ for Artinian rings
\cite[Theorems~3.7, 3.8]{MR15} and $2\le\gamma(\G(\Z_n))\le k$ for $n$ with $k$ prime divisors
\cite[Corollary~3.10]{MR15} become equalities when there are at least three maximal ideals. For a product of two local rings,
\cite[Proposition~3.3(c)]{MR15} gives $\gamma\in\{1,2\}$; Theorem~\ref{thm:gamma-k2} decides
between the two values for every ring with two maximal ideals. In \cite[Example~3.11(b)]{MR15} the
primes were shown to dominate $\G(R)$ for a principal ideal domain; by Theorem~\ref{thm:infinite},
no ring with infinitely many maximal ideals has a finite dominating set. The domatic number, the
independent domination number and the graph polynomials were not considered in \cite{MR15}.

The blow-up structure is due to \cite{MP11}. Gadge, Khandekar and Joshi \cite{GKJ24} realized
$\G(R)$, for Artinian $R$, as the zero-divisor graph of a lattice obtained from $2^k$ by replacing
its elements with chains, and blow-ups of Boolean graphs were used for comaximal ideal graphs in
\cite{YWLY14,GWY15}. For squarefree $n$, Rather \cite{Ra26d} described $\G(\Z_n)$ as a weighted
blow-up of $\TG_k$, with class sizes $\prod_{i\notin A}(p_i-1)$, and determined its vertex
connectivity. The class sizes \eqref{eq:weights-intro} for arbitrary rings, the isomorphism criterion
(Theorem~\ref{thm:iso}, Corollary~\ref{cor:iso-finite}) and the domination results are, to the best of
our knowledge, new. The papers \cite{Ra26a,Ra26b} also study zeros, unimodality and log-concavity of
the polynomials, which we do not consider.

Graph terminology follows \cite{HHS98,HY13}. A set $D$ of vertices is \emph{dominating} if every
vertex outside $D$ has a neighbour in $D$, \emph{total dominating} if every vertex has a neighbour
in $D$, and \emph{connected dominating} if it is dominating and induces a connected subgraph; the
minimum sizes are $\gamma,\gamma_t,\gamma_c$. The \emph{domatic number} $d$ is the maximum number of
pairwise disjoint dominating sets \cite{CH77}, and $i$ is the minimum size of an independent
dominating set \cite{GH13}.

\section{The type graph and the class sizes}\label{sec:type-graph}

\begin{lemma}\label{lem:adjacency}
An element $x\in R$ is a vertex of $\G(R)$ if and only if $\emptyset\neq\Lambda(x)\neq\Max(R)$, and
distinct vertices $x,y$ are adjacent if and only if $\Lambda(x)\cap\Lambda(y)=\emptyset$.
\end{lemma}

\begin{proof}
$x$ is a nonunit if and only if $Rx$ lies in a maximal ideal, that is, $\Lambda(x)\neq\emptyset$;
and $x\notin\Jac(R)$ if and only if $\Lambda(x)\neq\Max(R)$. Likewise $Rx+Ry\neq R$ if and only if
some maximal ideal contains $x$ and $y$.
\end{proof}

Hence $\G(R)=\emptyset$ exactly when $R$ is local; from now on $R$ is not local.

\begin{definition}\label{def:type-graph}
Let $k\ge2$, let $\Pk$ be the set of nonempty proper subsets of $[k]$, and write
$\co A=[k]\setminus A$. The \emph{type graph} $\TG_k$ has vertex set $\Pk$, with $A\sim B$ if
$A\cap B=\emptyset$. For a graph $H$ and nonzero cardinals $w(h)$, $h\in V(H)$, the \emph{blow-up}
$H[w]$ replaces each $h$ by an independent set $X_h$ of size $w(h)$, with $X_h$ and $X_{h'}$
completely joined when $hh'$ is an edge.
\end{definition}

The types $\{i\}$ are the \emph{singletons} and the types $\co{\{i\}}$ the \emph{co-singletons}. The
neighbours of $A$ in $\TG_k$ are the nonempty subsets of $\co A$. In particular $\co{\{i\}}$ is
adjacent only to $\{i\}$, the singletons form a $k$-clique, and $\TG_2=K_2$; see
Figure~\ref{fig:K3}.

\begin{figure}[ht]
\centering
\begin{tikzpicture}[every node/.style={draw,rounded corners,inner sep=2.5pt,font=\small,fill=white}]
  \node (1) at (90:1.1)  {$\{1\}$};
  \node (2) at (210:1.1) {$\{2\}$};
  \node (3) at (330:1.1) {$\{3\}$};
  \node (23) at (90:2.5)  {$\{2,3\}$};
  \node (13) at (210:2.5) {$\{1,3\}$};
  \node (12) at (330:2.5) {$\{1,2\}$};
  \draw (1)--(2)--(3)--(1);
  \draw (1)--(23); \draw (2)--(13); \draw (3)--(12);
\end{tikzpicture}
\caption{The type graph $\TG_3$.}
\label{fig:K3}
\end{figure}

Let $\Max(R)=\{M_1,\dots,M_k\}$ with $2\le k<\infty$, $F_i=R/M_i$ and $q_i=|F_i|$. For a vertex $x$
put $\tau(x)=\{i:x\in M_i\}\in\Pk$, and for $A\in\Pk$ let $V_A=\{x:\tau(x)=A\}$ and
$w_R(A)=|\Jac(R)|\prod_{i\notin A}|F_i^{\times}|$, a cardinal. Thus
\begin{equation}\label{eq:sing-cosing}
  w_R(\{i\})=|\Jac(R)|\prod_{j\neq i}|F_j^\times|,\qquad w_R(\co{\{i\}})=|\Jac(R)|\,|F_i^\times| .
\end{equation}

\begin{lemma}\label{lem:realization}
For every $A\in\Pk$, $V_A$ is a nonempty union of cosets of $\Jac(R)$, and $|V_A|=w_R(A)$.
\end{lemma}

\begin{proof}
By the Chinese Remainder Theorem, $\pi\colon R\to F_1\times\dots\times F_k$ is surjective with kernel
$\Jac(R)$. Since $x\in M_i$ if and only if $\pi(x)_i=0$, we have $V_A=\pi^{-1}(S_A)$, where $S_A$ is
the set of tuples vanishing exactly in the coordinates in $A$, and $|S_A|=\prod_{i\notin A}|F_i^\times|$.
\end{proof}

Nonemptiness of $V_A$ is \cite[Proposition~2.3]{MP11}. Note that $R$ is finite if and only if
$\Jac(R)$ and all $F_i$ are finite.

\begin{theorem}[cf.\ {\cite[Theorem~2.4]{MP11}}]\label{thm:structure}
If $R$ has exactly $k$ maximal ideals, $2\le k<\infty$, then $\G(R)\cong\TG_k[w_R]$, with $V_A$
mapped onto $X_A$.
\end{theorem}

\begin{proof}
This follows from Lemmas~\ref{lem:adjacency} and \ref{lem:realization}.
\end{proof}

For $k=2$ this gives $\G(R)=K_{a,b}$ \cite[Theorem~2.2]{MSSY08}. For $k\ge3$ a vertex of type
$\co{\{k\}}$ and a vertex of type $\co{\{1\}}$ are at distance $3$, which recovers
\cite[Proposition~3.3(b)]{MR15}.

\begin{lemma}\label{lem:twins}
Let $k\ge2$ and let $w$ take nonzero values.
\begin{enumerate}[label=\textup{(\roman*)}]
\item Two vertices of $\TG_k[w]$ have the same neighbourhood if and only if they lie in the same
  class.
\item $\Aut(\TG_k)\cong S_k$, acting naturally on $\Pk$.
\end{enumerate}
\end{lemma}

\begin{proof}
(i) If $i\in A\setminus B$, the vertices of $X_{\{i\}}$ are adjacent to those of $X_B$ but not to
those of $X_A$. (ii) The degree of $A$ is $2^{k-|A|}-1$, so an automorphism permutes the singletons,
say $\{i\}\mapsto\{\sigma(i)\}$. Since $A=\{i:\{i\}=A\text{ or }\{i\}\not\sim A\}$, it maps $A$ to
$\sigma(A)$.
\end{proof}

\begin{theorem}\label{thm:iso}
Let $R,R'$ have $k,k'$ maximal ideals, $2\le k,k'<\infty$. Then $\G(R)\cong\G(R')$ if and only if
$k=k'$ and $w_{R'}(\sigma(A))=w_R(A)$ for all $A\in\Pk$, for some $\sigma\in S_k$ acting on $\Pk$
by $A\mapsto\sigma(A)$.
\end{theorem}

\begin{proof}
By Lemma~\ref{lem:twins}(i) an isomorphism maps classes onto classes, so it induces an isomorphism
$\TG_k\to\TG_{k'}$ preserving class sizes. Hence $k=k'$, and Lemma~\ref{lem:twins}(ii) applies.
\end{proof}

\begin{corollary}\label{cor:iso-finite}
Let $R,R'$ be finite rings, each with at least three maximal ideals. Then $\G(R)\cong\G(R')$ if and only if
$|\Jac(R)|=|\Jac(R')|$ and the residue fields have the same orders, counted with multiplicity. In
particular, if $n$ and $n'$ have at least three prime divisors, then $\G(\Z_n)\cong\G(\Z_{n'})$ if
and only if $n=n'$.
\end{corollary}

\begin{proof}
Let $j=|\Jac(R)|$ and $Q=\prod_i(q_i-1)$. By \eqref{eq:sing-cosing},
$w_R(\{1\})w_R(\co{\{1\}})=j^2Q$ and $\prod_iw_R(\co{\{i\}})=j^kQ$, so $w_R$ determines $j^{k-2}$,
hence $j$, and then $q_i=1+w_R(\co{\{i\}})/j$. Conversely these data determine $w_R$ up to
permutation. For $\Z_n$ the $q_i$ are the prime divisors $p_i$ of $n$ and $j=n/\rad(n)$, so
$n=j\prod_ip_i$ is determined.
\end{proof}

This fails for $k=2$: $\G(\Z_4\times\Z_4)\cong\G(\FF_5\times\FF_5)\cong K_{4,4}$. For infinite rings
the graph retains even less. Moconja and Petrovi\'c \cite[Example~4.1]{MP11} observed that the
graphs of the semilocal rings $\Z_{(p,q)}$ are all isomorphic; by \eqref{eq:weights-intro} all
their classes have size $\aleph_0$, and by Theorem~\ref{thm:iso} the class sizes are the only
invariant.

\section{Domination in weighted type graphs}\label{sec:domination}

In this section $G=\TG_k[w]$ with $k\ge2$ and arbitrary nonzero cardinals $w(A)$; a vertex in $X_A$
has type $A$. Put $Y_i=X_{\{i\}}\cup X_{\co{\{i\}}}$. A \emph{transversal} is a set containing exactly
one vertex of each $X_{\{i\}}$. We use three facts about $\TG_k$:
\begin{enumerate}[label=(F\arabic*)]
\item\label{F1} each type $A$ is adjacent to $\{i\}$ for every $i\notin A$;
\item\label{F2} for $k\ge3$ the $2k$ types $\{i\}$, $\co{\{i\}}$ are distinct, so $Y_1,\dots,Y_k$ are
  disjoint; for $k=2$, $\co{\{1\}}=\{2\}$ and $Y_1=Y_2$;
\item\label{F3} for $k=3$ every type is a singleton or a co-singleton; for $k\ge4$ the types $A$
  with $2\le|A|\le k-2$ are neither.
\end{enumerate}

\begin{lemma}\label{lem:key}
Let $D$ be a dominating set. If $D\cap X_{\{i\}}=\emptyset$, then $X_{\co{\{i\}}}\subseteq D$. In
particular $D$ meets every $Y_i$.
\end{lemma}

\begin{proof}
Every neighbour of a vertex of type $\co{\{i\}}$ has type $\{i\}$.
\end{proof}

\begin{theorem}\label{thm:gamma-k3}
If $k\ge3$, then $\gamma(G)=\gamma_t(G)=\gamma_c(G)=\omega(G)=\chi(G)=k$, and every transversal is a
dominating clique.
\end{theorem}

For $\G(R)$ the equalities $\omega=\chi=k$ are in \cite{MP11}; we include the short argument.

\begin{proof}
By \ref{F2} and Lemma~\ref{lem:key}, $\gamma\ge k$, hence also $\gamma_t,\gamma_c\ge k$. A transversal is a clique of size $k$, and it
dominates by \ref{F1}; being a clique, it is also total and connected dominating, and $\omega\ge k$. The types in a
clique are pairwise disjoint, so $\omega\le k$, and colouring each vertex by the least element of
its type is proper, so $\chi\le k$.
\end{proof}

\begin{proposition}\label{prop:k2-abstract}
Let $k=2$, so $G=K_{a,b}$ with $a=w(\{1\})$ and $b=w(\{2\})$. Then $\gamma=\gamma_c=1$ if
$\min\{a,b\}=1$ and $\gamma=\gamma_c=2$ otherwise, $\gamma_t=2$, $i=\min\{a,b\}$, and
$d=\max\{\min\{a,b\},2\}$ if $\min\{a,b\}$ is finite, and $d=\min\{a,b\}$ otherwise.
\end{proposition}

\begin{proof}
Only $d$ needs comment. A set dominates $K_{a,b}$ if and only if it meets both parts or contains a
part. If one of several disjoint dominating sets contains a part, every other one contains the other
part, so there are at most two; otherwise there are at most $\min\{a,b\}$. Both bounds are attained.
\end{proof}

\subsection{Minimum dominating sets}

\begin{theorem}\label{thm:gamma-sets}
Let $k\ge3$. A set $D$ is a minimum dominating set if and only if $D=\{v_1,\dots,v_k\}$ with
$v_i\in Y_i$ and $E=\{i:v_i\in X_{\co{\{i\}}}\}$ satisfies
\begin{enumerate}[label=\textup{(\alph*)}]
\item $w(\co{\{i\}})=1$ for all $i\in E$;
\item $|E|\le1$ if $k\ge4$;
\item $w(\{j\})=1$ for the index $j\notin E$, if $k=3$ and $|E|=2$.
\end{enumerate}
If all $w(\{i\})$ are finite, the number of minimum dominating sets is
$\sum_E\prod_{i\notin E}w(\{i\})$, the sum running over the sets $E\subseteq[k]$ satisfying
\textup{(a)--(c)}.
\end{theorem}

\begin{proof}
By Lemma~\ref{lem:key} and \ref{F2}, a dominating set of size $k$ meets each $Y_i$ in exactly one
vertex $v_i$ and has no other elements.

\emph{Necessity.} (a) If $i\in E$, then $X_{\co{\{i\}}}\subseteq D$ by Lemma~\ref{lem:key}, so
$X_{\co{\{i\}}}=\{v_i\}$. (b) Let $k\ge4$ and $i\ne l$ in $E$. Then $B=[k]\setminus\{i,l\}$ is neither
a singleton nor a co-singleton, so $X_B\cap D=\emptyset$. The neighbours of $X_B$ have types
$\{i\}$, $\{l\}$ or $\{i,l\}$, and none of these occurs in $D$, since $\{i,l\}$ is not a co-singleton.
(c) Let $k=3$, $E=\{i,l\}$ and $j\notin E$. The neighbours of $X_{\{j\}}$ have types $\{i\}$, $\{l\}$
or $\{i,l\}=\co{\{j\}}$, none of which occurs in $D$; so $X_{\{j\}}=\{v_j\}$.

\emph{Sufficiency.} Let $y\notin D$ have type $A$; by \ref{F3} there are four cases. If $A=\{j\}$
with $j\in E$, then $y\sim v_j$. If $A=\{j\}$ with $j\notin E$, any $i\notin E\cup\{j\}$ gives
$y\sim v_i$; such an $i$ exists unless $|E|=k-1$, which by (b) means $k=3$ and $|E|=2$, and then (c)
gives $y\in D$. If $A=\co{\{i\}}$, then $i\notin E$ by (a), and $y\sim v_i$. If $2\le|A|\le k-2$, then
$k\ge4$ and $|E|\le1<2\le|[k]\setminus A|$, so there is $j\in[k]\setminus(A\cup E)$; then $v_j$ has
type $\{j\}$, which is disjoint from $A$, and $y\sim v_j$.

The count follows, since $v_i$ is forced for $i\in E$ and ranges over $X_{\{i\}}$ for $i\notin E$.
\end{proof}

For example, every class of $\G(\Bool3)$ has one vertex and all $8$ sets $E$ are admissible, while
$\G(\Bool4)$ has $1+4=5$ minimum dominating sets.

\subsection{The domatic number}

\begin{theorem}\label{thm:domatic}
If $k\ge3$ and $m=\min_iw(\{i\})$ is finite, then $d(G)=m+1$.
\end{theorem}

\begin{proof}
Let $w(\{i_0\})=m$. By Lemma~\ref{lem:key}, each member of a family of disjoint dominating sets
meets $X_{\{i_0\}}$ or contains $X_{\co{\{i_0\}}}$; at most $m$ members do the former and at most
one the latter. Conversely, take $m$ disjoint transversals $S_1,\dots,S_m$ and let
$S_{m+1}=V(G)\setminus(S_1\cup\dots\cup S_m)$. The vertices outside $S_{m+1}$ have singleton
types. By \ref{F2}, $S_{m+1}$ contains every co-singleton class, and a vertex of type $\{i\}$ is
adjacent to $X_{\co{\{i\}}}$; so $S_{m+1}$ is dominating.
\end{proof}

If $m$ is infinite, the same argument gives $d(G)=m$. For $k=2$ the answer is different
(Proposition~\ref{prop:k2-abstract}), because $Y_1=Y_2$.

\section{Independent domination}\label{sec:independent}

We keep $G=\TG_k[w]$ with $k\ge2$. For $\F\subseteq\Pk$ put $X_\F=\bigcup_{A\in\F}X_A$ and
$w(\F)=\sum_{A\in\F}w(A)$. A family $\F\subseteq\Pk$ is \emph{intersecting} if any two members
meet, and \emph{maximal intersecting} if it is maximal among intersecting subfamilies of $\Pk$.
For arbitrary positive real weights $w$ we read $i(\TG_k[w])$ as the minimum of $w(\F)$ over the
maximal intersecting families $\F$; by Theorem~\ref{thm:indep-char} below this agrees with the
graph-theoretic definition when the weights are cardinals.

\begin{theorem}\label{thm:indep-char}
The independent dominating sets of $G$ are the sets $X_\F$ with $\F\subseteq\Pk$ maximal
intersecting. Hence $i(G)$ is the minimum of $w(\F)$ over these families.
\end{theorem}

\begin{proof}
Let $D$ be independent and dominating, with set of types $\F$; then $\F$ is intersecting. If $B\in\Pk$
meets every member of $\F$, no vertex of $X_B$ has a neighbour in $D$, so $X_B\subseteq D$. Hence
$\F$ is maximal, and $D=X_\F$. Conversely, if $\F$ is maximal intersecting, then $X_\F$ is
independent, and a vertex of type $B\notin\F$ is adjacent to $X_A$ for some $A\in\F$ disjoint from
$B$.
\end{proof}

\begin{lemma}\label{lem:one-of-each}
A family $\F\subseteq\Pk$ is maximal intersecting if and only if it is intersecting and contains
one of $A,\co A$ for every $A\in\Pk$. In particular $|\F|=2^{k-1}-1$.
\end{lemma}

\begin{proof}
If $\F$ is maximal intersecting and contains neither $A$ nor $\co A$, it contains some
$B\subseteq\co A$ and some $C\subseteq A$, which are disjoint. The converse holds because
$C\cap\co C=\emptyset$.
\end{proof}

\begin{corollary}\label{cor:i-bounds}
With the sum over the $2^{k-1}-1$ complementary pairs $\{A,\co A\}$,
\[
  \sum_{\{A,\co A\}}\min\{w(A),w(\co A)\}\ \le\ i(G)\ \le\ \sum_{1\in A\in\Pk}w(A).
\]
\end{corollary}

\begin{proof}
The lower bound follows from Lemma~\ref{lem:one-of-each}, the upper bound from the family
$\{A:1\in A\}$.
\end{proof}

\begin{theorem}\label{thm:i-closed}
If $w(A)=c\prod_{j\notin A}a_j$ for real numbers $c>0$ and $a_1,\dots,a_k\ge1$, then the lower bound in
Corollary~\ref{cor:i-bounds} is attained:
\[
 i(G)=\frac c2\sum_{A\in\Pk}\min\Bigl\{\prod_{j\in A}a_j,\ \prod_{j\notin A}a_j\Bigr\}.
\]
\end{theorem}

\begin{proof}
Let $s_i=\ln a_i\ge0$ and $s(X)=\sum_{i\in X}s_i$ for $X\subseteq[k]$. Since
$\ln w(A)=\ln c+s(\co A)$, we have $w(A)\le w(\co A)$ if and only if $s(A)\ge s(\co A)$. Let
$\Delta$ be the least nonzero value of $|s(A)-s(\co A)|$, $A\in\Pk$, or $\Delta=1$ if there is
none. Choose $\delta>0$ with $\delta2^{k+1}<\Delta$, and put $\varepsilon(X)=\delta\sum_{i\in X}2^i$
and $\phi=s+\varepsilon$. Thus $\phi$ orders the sets as $s$ does, ties being broken consistently
by $\varepsilon$.

For every $A\in\Pk$ we have $0<|\varepsilon(A)-\varepsilon(\co A)|<\Delta$: the difference is
nonzero because distinct sets have distinct binary sums, and it is at most
$\delta\sum_{i=1}^k2^i<\delta2^{k+1}$. If $s(A)=s(\co A)$, then
$\phi(A)-\phi(\co A)=\varepsilon(A)-\varepsilon(\co A)\neq0$. Otherwise
$|s(A)-s(\co A)|\ge\Delta$, so $\phi(A)-\phi(\co A)$ is nonzero and has the sign of $s(A)-s(\co A)$.
In both cases $\phi(A)\ne\phi(\co A)$, and $\phi(A)>\phi(\co A)$ implies $s(A)\ge s(\co A)$, that is,
$w(A)\le w(\co A)$. Hence $\F=\{A\in\Pk:\phi(A)>\phi(\co A)\}$ contains exactly one set of each pair
$\{A,\co A\}$, namely one of smaller weight.

Since $s_i\ge0$ and $\delta>0$, $\phi$ is strictly increasing under inclusion. If $A,B\in\F$ were
disjoint, then $B\subseteq\co A$ and $A\subseteq\co B$, so
\[
  \phi(B)\le\phi(\co A)<\phi(A)\le\phi(\co B)<\phi(B),
\]
a contradiction. So $\F$ is intersecting, hence maximal by Lemma~\ref{lem:one-of-each}, and
$w(\F)$ equals the lower bound of Corollary~\ref{cor:i-bounds}.
\end{proof}

\begin{remark}\label{rem:i-fails}
The hypothesis is needed. For $k=3$ and $w(A)=2^{|A|}$, the maximal intersecting families are
$\{\{1,2\},\{1,3\},\{2,3\}\}$, of weight $12$, and the three families $\{\{i\},\{i,j\},\{i,l\}\}$, of
weight $10$. So $i(G)=10$, while the lower bound is $6$. These weights occur in
Example~\ref{ex:divisor-lattice}.
\end{remark}

The maximal intersecting families on $[k]$ correspond to the self-dual monotone Boolean functions
of $k$ variables, and their number is known for $k\le9$ \cite{BMMV13}. For completeness we give a
short argument that it is at least doubly exponential in $k$. Let $k=2\ell+1\ge3$, let $\mathcal B$ be
the set of $(\ell+1)$-subsets of $[k]$ containing $\{1,2\}$, so $|\mathcal B|=\binom{2\ell-1}{\ell-1}$, and
for $\mathcal S\subseteq\mathcal B$ let $\F_{\mathcal S}$ consist of the sets in $\Pk$ of size at
least $\ell+1$ not in $\mathcal S$, together with the sets $[k]\setminus B$ for $B\in\mathcal S$. Two sets of size at
least $\ell+1$ meet. A set $C$ is disjoint from $[k]\setminus B$ only if $C\subseteq B$; for
$|C|\ge \ell+1$ this means $C=B\notin\F_{\mathcal S}$. For $B,B'\in\mathcal S$ we have
$|B\cup B'|\le 2\ell<k$, so $([k]\setminus B)\cap([k]\setminus B')\neq\emptyset$. Thus
$\F_{\mathcal S}$ is intersecting and contains one set of each complementary pair, so it is maximal
by Lemma~\ref{lem:one-of-each}. Since $B\in\F_{\mathcal S}$ if and only if $B\notin\mathcal S$, the
$2^{|\mathcal B|}$ families $\F_{\mathcal S}$ are distinct, and there are at least
$2^{\binom{2\ell-1}{\ell-1}}$ maximal intersecting families on $[2\ell+1]$. This extends to even $k$, since
$\F\mapsto\{A\in\mathcal P_{k+1}:A\cap[k]\in\F\cup\{[k]\}\}$ maps the maximal intersecting families
on $[k]$ injectively to those on $[k+1]$.

\section{Consequences for rings}\label{sec:rings}

By Theorem~\ref{thm:structure}, the results of Sections~\ref{sec:domination} and
\ref{sec:independent} apply to $\G(R)$ with $w=w_R$.

\subsection{Finitely many maximal ideals}

\begin{theorem}\label{thm:gamma-k2}
If $R$ has two maximal ideals, then $\G(R)=K_{a,b}$ with $a=w_R(\{1\})$ and $b=w_R(\{2\})$, and
$\gamma(\G(R))=1$ if $R\cong\FF_2\times F$ for a field $F$, and $\gamma(\G(R))=2$ otherwise.
\end{theorem}

\begin{proof}
By \eqref{eq:sing-cosing}, $a=1$ if and only if $\Jac(R)=0$ and $R/M_2\cong\FF_2$, and then
$R\cong R/M_1\times\FF_2$ by the Chinese Remainder Theorem; similarly for $b$. Apply Proposition~\ref{prop:k2-abstract}.
\end{proof}

\begin{corollary}\label{cor:rings-finite-max}
Let $R$ have $k\ge3$ maximal ideals.
\begin{enumerate}[label=\textup{(\roman*)}]
\item $\gamma=\gamma_t=\gamma_c=k$ for $\G(R)$ (and $\omega=\chi=k$ \cite{MP11}), and elements
  $x_i\in M_i\setminus\bigcup_{j\ne i}M_j$, $1\le i\le k$, form a minimum dominating clique.
\item In Theorem~\ref{thm:gamma-sets}, condition (a) holds for $i$ if and only if $\Jac(R)=0$ and
  $R/M_i\cong\FF_2$; condition (c) then holds automatically.
\item If $m=\min_iw_R(\{i\})$ is finite, then $d(\G(R))=m+1$.
\end{enumerate}
\end{corollary}

\begin{proof}
Apply Theorems~\ref{thm:gamma-k3}, \ref{thm:gamma-sets} and \ref{thm:domatic}, using
\eqref{eq:sing-cosing}. In (ii), if $|E|=2$ then $\Jac(R)=0$ and two residue fields are $\FF_2$, so
$w_R(\{j\})=1$.
\end{proof}

\begin{corollary}\label{cor:i-rings}
Let $R$ be finite and non-local, with residue fields of orders $q_1,\dots,q_k$. Then
\[
  i(\G(R))=\frac{|\Jac(R)|}{2}\sum_{A\in\Pk}\min\Bigl\{\prod_{j\in A}(q_j-1),\ \prod_{j\notin A}(q_j-1)\Bigr\}.
\]
If all $q_i$ equal $q$, then
\[
  i(\G(R))=|\Jac(R)|\Bigl(\sum_{k/2<s<k}\tbinom ks(q-1)^{k-s}+\varepsilon_k\tfrac12\tbinom{k}{k/2}(q-1)^{k/2}\Bigr),
\]
where $\varepsilon_k$ is $1$ for even $k$ and $0$ for odd $k$. In particular
$i(\G(\Bool k))=2^{k-1}-1$, which exceeds $\gamma(\G(\Bool k))=k$ for $k\ge4$.
\end{corollary}

\begin{proof}
Apply Theorem~\ref{thm:i-closed} with $c=|\Jac(R)|$ and $a_i=q_i-1$. If all $q_i=q$, the lighter set
in a pair $\{A,\co A\}$ is the larger one, of weight $|\Jac(R)|(q-1)^{k-|A|}$.
\end{proof}

\begin{corollary}\label{cor:Zn}
Let $n=p_1^{t_1}\cdots p_k^{t_k}$ with $k\ge2$ distinct primes, $\rad(n)=p_1\cdots p_k$ and
$\bar n=n/\rad(n)$. The class of $A$ in $\G(\Z_n)$ has $\bar n\,\varphi(\prod_{i\notin A}p_i)$
elements, and:
\begin{enumerate}[label=\textup{(\roman*)}]
\item $\gamma(\G(\Z_n))=k$ if $k\ge3$; for $k=2$, $\gamma(\G(\Z_n))=1$ if and only if $n=2p$ with $p$
  an odd prime;
\item $d(\G(\Z_n))=1+\varphi(n)/(P-1)$ if $k\ge3$, where $P=\max_ip_i$;
\item $i(\G(\Z_n))=\dfrac{\bar n}{2}\displaystyle\sum_{e\mid\rad(n),\ 1<e<\rad(n)}
  \min\bigl\{\varphi(e),\varphi(\rad(n)/e)\bigr\}$.
\end{enumerate}
\end{corollary}

\begin{proof}
Here $|\Jac(\Z_n)|=\bar n$, $q_i=p_i$ and $w(\{i\})=\varphi(n)/(p_i-1)$.
\end{proof}

For example, $\G(\Z_{210})$ has $161$ vertices, and $\gamma=4$, $d=9$ and $i=25$.

\subsection{Graph polynomials}
Let $G=\TG_k[w]$ be as in Sections~\ref{sec:domination} and~\ref{sec:independent}. For a finite
graph $H$ let $I(H;x)$, $D_i(H;x)$ and $D(H;x)$ be the generating functions, by size, of
the independent sets (including $\emptyset$), the independent dominating sets and the dominating
sets. For $\mathcal T\subseteq\Pk$ let $U(\mathcal T)$ be the set of $A\in\Pk$ that meet every
member of $\mathcal T$.

\begin{theorem}\label{thm:polys}
If $w$ is finite-valued, then
\begin{align*}
I(G;x)   &= \sum_{\F\text{ intersecting}}\ \prod_{A\in\F}\bigl((1+x)^{w(A)}-1\bigr),\qquad
D_i(G;x) = \sum_{\F\text{ maximal intersecting}} x^{w(\F)},\\
D(G;x)   &= \sum_{U(\mathcal T)\subseteq\mathcal T\subseteq\Pk}\
            \prod_{A\in U(\mathcal T)}x^{w(A)}\prod_{A\in\mathcal T\setminus U(\mathcal T)}\bigl((1+x)^{w(A)}-1\bigr).
\end{align*}
\end{theorem}

\begin{proof}
Classify a vertex set $D$ by $\mathcal T=\{A:D\cap X_A\ne\emptyset\}$. It is independent if and only
if $\mathcal T$ is intersecting; this gives $I(G;x)$, since $D$ meets each $X_A$, $A\in\mathcal T$,
in a nonempty set. The types in $U(\mathcal T)$ are exactly those having no neighbour in $D$, so $D$
is dominating if and only if $X_A\subseteq D$ for every $A\in U(\mathcal T)$; in particular
$U(\mathcal T)\subseteq\mathcal T$. For $A\in\mathcal T\setminus U(\mathcal T)$ the vertices of
$X_A$ are dominated anyway, so $D\cap X_A$ may be any nonempty subset of $X_A$. This gives
$D(G;x)$. The formula for $D_i$ is Theorem~\ref{thm:indep-char}.
\end{proof}

The sum for $D_i$ has one term for each maximal intersecting family, and terms of equal weight
combine; for $\Bool k$ all weights are $1$, and $D_i$ is a single monomial.

\begin{corollary}\label{cor:full-graph-polys}
Let $R$ be finite and non-local, with $|R|=N$, $|U(R)|=u$ and $|\Jac(R)|=j$. Then
\begin{align*}
I(\Gamma(R);x)   &= ux + (1+x)^{j}\,I(\G(R);x),\qquad
D_i(\Gamma(R);x) = ux + x^{j}\,D_i(\G(R);x),\\
D(\Gamma(R);x)   &= \bigl((1+x)^{u}-1\bigr)(1+x)^{N-u} + x^{j}\,D(\G(R);x).
\end{align*}
\end{corollary}

\begin{proof}
A unit is adjacent to every vertex, and $\Jac(R)$ is an independent set adjacent only to units. So
an independent set containing a unit is a singleton, every set containing a unit is dominating, and
a dominating set without units contains $\Jac(R)$ and meets $\G(R)$ in a dominating set of $\G(R)$.
\end{proof}

\begin{example}
Let $n=30$ and index the types by subsets of $\{2,3,5\}$. The classes of
$\{2\}$, $\{3\}$, $\{5\}$, $\{2,3\}$, $\{2,5\}$, $\{3,5\}$ have sizes $8,4,2,4,2,1$. The maximal intersecting
families are the family of the three $2$-sets, of weight $7$, and the families
$\{\{s\},\{s,t\},\{s,u\}\}$ with $\{s,t,u\}=\{2,3,5\}$, of weights $14$, $9$, $5$ for $s=2,3,5$.
By Corollary~\ref{cor:full-graph-polys}, with $|U(\Z_{30})|=8$ and $|\Jac(\Z_{30})|=1$, the
exponents in $\Gamma(\Z_{30})$ exceed these weights by $1$:
$D_i(\Gamma(\Z_{30});x)=8x+x^6+x^8+x^{10}+x^{15}$, as in \cite{Ra26a}. For $n=p_1p_2p_3$ one obtains
in the same way the exponents $p_1+p_2+p_3-2$, $p_1p_2$, $p_1p_3$, $p_2p_3$ found in \cite{Ra26a}.
\end{example}

\subsection{Infinitely many maximal ideals}

\begin{theorem}\label{thm:infinite}
If $\Max(R)$ is infinite, then $\G(R)$ has no finite dominating set.
\end{theorem}

\begin{proof}
Let $D=\{z_1,\dots,z_n\}$ be a nonempty set of vertices. Choose maximal ideals $N_j\ni z_j$, and
distinct maximal ideals $P_1,\dots,P_{n+1}$ not among the $N_j$. For $1\le t\le n+1$ put
\[
  I_t=N_1N_2\cdots N_n\cdot\prod_{s\neq t}P_s .
\]
Distinct maximal ideals are not contained in one another, so no factor of $I_t$ is contained in
$P_t$; since $P_t$ is prime, neither is the product $I_t$. Choose $y_t\in I_t\setminus P_t$. Each $y_t$ lies in every $N_j$ and not in $P_t$, so it is a
vertex; and the $y_t$ are distinct,
since $y_u\in P_t$ for $u\ne t$. Some $y_t$ is not in $D$, and it shares the maximal ideal $N_j$ with
each $z_j$.
\end{proof}

Together with Theorems~\ref{thm:gamma-k3} and \ref{thm:gamma-k2}, this proves Theorem~A.

\begin{example}\label{ex:boolean-infinite}
(i) If every maximal ideal $M$ contains an element $x_M$ lying in no other maximal ideal, then
$\{x_M:M\in\Max(R)\}$ dominates $\G(R)$, since a vertex $y\notin M$ is adjacent to $x_M$.
For $\Z$, or any countable PID with infinitely many primes, this and Theorem~\ref{thm:infinite} give
$\gamma(\G(R))=\aleph_0=|\Max(R)|$, recovering \cite[Example~3.11(b)]{MR15}.

(ii) Let $B$ be the free Boolean ring on countably many generators. It is countable, and its maximal
ideals are the kernels of the homomorphisms $B\to\FF_2$, which correspond to the assignments of the
generators; so $|\Max(B)|=2^{\aleph_0}$. Since the vertex set is countable,
Theorem~\ref{thm:infinite} gives $\gamma(\G(B))=\aleph_0<|\Max(B)|$. The same argument gives $\gamma(\G(R))=\aleph_0$ for every
countable $R$ with infinitely many maximal ideals.
\end{example}

\section{Comaximal graphs of lattices}\label{sec:lattices}

Let $L$ be a bounded lattice. An ideal is a nonempty down-set closed under finite joins; by Zorn's
lemma every proper ideal lies in a maximal one. Put $\Jac(L)=\bigcap\Max(L)$. The comaximal graph of
$L$ \cite{AK14} has $a\neq b$ adjacent if $a\vee b=1$, and as for rings we consider its subgraph
$\G(L)$ on $L\setminus(\{1\}\cup\Jac(L))$. For finite $L$, Afkhami and Khashyarmanesh proved
$\chi(\G(L))=\omega(\G(L))=|C(L)|$, where $C(L)$ is the set of coatoms \cite[Theorem~3.4]{AK14};
see also \cite{AKS15}. Domination in $\G(L)$ does not seem to have been studied.

Since $a\vee b\neq1$ if and only if some maximal ideal contains $a$ and $b$,
Lemma~\ref{lem:adjacency} holds for $\G(L)$. We assume that every maximal ideal of $L$ is prime,
that is, $u\wedge v\in M$ implies $u\in M$ or $v\in M$. This holds for distributive lattices
\cite{Gr11} and also, for instance, for the pentagon $N_5$.

\begin{lemma}\label{lem:lattice-realization}
If the maximal ideals $M_1,\dots,M_k$ of $L$ are prime, every $A\in\Pk$ is the type of some
$x\in L$.
\end{lemma}

\begin{proof}
For $i\in A$ and $j\notin A$ choose $u_{ij}\in M_i\setminus M_j$, and put
$y_i=\bigvee_{j\notin A}u_{ij}$ and $x=\bigwedge_{i\in A}y_i$. Then $y_i\in M_i$, and $y_i\notin M_j$
for $j\notin A$ because $y_i\ge u_{ij}$. Hence $x\in M_i$ for $i\in A$, and $x\notin M_j$ for
$j\notin A$ because $M_j$ is prime.
\end{proof}

\begin{theorem}\label{thm:lattice}
Let $L$ be a bounded lattice whose maximal ideals are prime.
\begin{enumerate}[label=\textup{(\roman*)}]
\item If $2\le|\Max(L)|=k<\infty$, then $\G(L)\cong\TG_k[w_L]$, where $w_L(A)$ is the number of
  elements of type $A$. Hence Theorems~\ref{thm:gamma-k3}, \ref{thm:gamma-sets}, \ref{thm:domatic}
  and \ref{thm:indep-char} and Proposition~\ref{prop:k2-abstract} apply to $\G(L)$, and so does
  Theorem~\ref{thm:polys} if all $w_L(A)$ are finite.
\item If $\Max(L)$ is infinite, then $\G(L)$ has no finite dominating set.
\end{enumerate}
\end{theorem}

\begin{proof}
(i) follows from Lemma~\ref{lem:lattice-realization}. For (ii), let $D=\{z_1,\dots,z_n\}$ be a
nonempty set of vertices, choose maximal ideals $N_j\ni z_j$ and distinct maximal ideals
$P_1,\dots,P_{n+1}$ not among the $N_j$, and argue as in Theorem~\ref{thm:infinite} with
$y_t=\bigwedge_jv_j\wedge\bigwedge_{s\ne t}v'_s$, where $v_j\in N_j\setminus P_t$ and
$v'_s\in P_s\setminus P_t$. Then $y_t\le v_j$ and $y_t\le v'_s$, so $y_t$ lies in every $N_j$ and every $P_s$ with $s\neq t$,
since ideals are down-sets. None of the elements $v_j,v'_s$ lies in $P_t$, so $y_t\notin P_t$,
since $P_t$ is prime.
\end{proof}

\begin{remark}\label{rem:lattice-credit}
Primeness cannot be dropped. In the diamond $M_3$ with atoms $a,b,c$, the maximal ideals
${\downarrow}a,{\downarrow}b,{\downarrow}c$ are not prime, since $b\wedge c=0\in{\downarrow}a$. No
element lies in exactly two of them, so $\G(M_3)$ is a triangle and $\gamma(\G(M_3))=1<3$.
\end{remark}

In a divisor lattice the class sizes do not decrease with the type, and Theorem~\ref{thm:i-closed}
does not apply in general (Remark~\ref{rem:i-fails}). For such weights the minimum is attained by a
family $\{A\in\Pk:i\in A\}$ with $t_i$ minimal. For $r\in[0,1]^k$ and a family $\mathcal A$ of
subsets of $[k]$ put $\mu_r(\mathcal A)=\sum_{X\in\mathcal A}\prod_{i\in X}r_i\prod_{i\notin X}(1-r_i)$.

\begin{theorem}[\cite{FFFLO86,Fr08}]\label{thm:product}
If $0<r\le\frac12$ and $\mathcal A$ is an intersecting family of subsets of $[k]$, then
$\mu_{(r,\dots,r)}(\mathcal A)\le r$.
\end{theorem}

\begin{proposition}\label{prop:i-increasing}
Let $w(A)=c\prod_{j\in A}t_j$ with $c>0$ and $t_1,\dots,t_k\ge1$, and let $t_1=\min_jt_j$. Then
$\{A\in\Pk:1\in A\}$ has the smallest weight among the maximal intersecting families, and
\[
  i(\TG_k[w])=c\Bigl(t_1\prod_{j=2}^k(1+t_j)-\prod_{j=1}^kt_j\Bigr).
\]
\end{proposition}

\begin{proof}
Put $p_i=t_i/(1+t_i)$, $r_i=1-p_i\le\frac12$ and $T=\prod_j(1+t_j)$. Then $w(A)=cT\mu_p(\{A\})$ and
$\mu_p(\{\co X\})=\mu_r(\{X\})$ for $X\subseteq[k]$. Let $\F$ be maximal intersecting and
$\F'=\F\cup\{[k]\}$. By Lemma~\ref{lem:one-of-each}, $\F'$ is intersecting and contains exactly one
set of each pair $\{X,\co X\}$ of subsets of $[k]$, so $\mu_p(\F')+\mu_r(\F')=1$ and
\[
  w(\F)=cT\bigl(\mu_p(\F')-\mu_p(\{[k]\})\bigr)=cT\bigl(1-\mu_r(\F')\bigr)-c\prod_jt_j .
\]
If $A\in\F'$ and $A\subseteq B$, then $\co B$ is disjoint from $A$, so $\co B\notin\F'$ and
$B\in\F'$; thus
$\F'$ is an up-set, and $\mu_r(\F')$ is non-decreasing in each $r_i$. Hence, with $r_1=\max_ir_i$,
$\mu_r(\F')\le\mu_{(r_1,\dots,r_1)}(\F')\le r_1$ by Theorem~\ref{thm:product}, and $w(\F)\ge cT(1-r_1)-c\prod_jt_j$, which is the
stated value. The family $\{A:1\in A\}$ attains it, since $\mu_r(\{X:1\in X\})=r_1$.
\end{proof}

\begin{example}\label{ex:divisor-lattice}
Let $L$ be the lattice of divisors of $n=p_1^{t_1}\cdots p_k^{t_k}$. Its maximal ideals are the sets
of divisors of $n/p_i$, a divisor $e$ has type $\{i:v_{p_i}(e)<t_i\}$, and
$w_L(A)=\prod_{i\in A}t_i$, which is non-decreasing in $A$. For $k\ge3$, $\gamma(\G(L))=k$ and
$d(\G(L))=1+\min_jt_j$. If $t_1=\min_jt_j$, Proposition~\ref{prop:i-increasing} gives
$i(\G(L))=t_1\prod_{j\ge2}(1+t_j)-\prod_jt_j$; for $n=900$ this is $10$, above the lower bound $6$
of Corollary~\ref{cor:i-bounds}. If all $t_i=1$, then $\G(L)\cong\G(\Bool k)$.
\end{example}

\section{Open problems}\label{sec:problems}

\begin{problem}
Characterize the rings with infinitely many maximal ideals for which $\gamma(\G(R))=|\Max(R)|$. Is
$\gamma(\G(R))\le|\Max(R)|$ always?
\end{problem}

\begin{problem}
Determine $i(\TG_k[w])$ for $w(A)=c\prod_{j\in A}t_j$ with $c,t_1,\dots,t_k>0$, when some $t_j$ are
smaller and some larger than $1$. The cases where all $t_j\le1$ or all $t_j\ge1$ are covered by
Theorem~\ref{thm:i-closed} (with $c\prod_jt_j$ in place of $c$ and $a_j=1/t_j$) and by
Proposition~\ref{prop:i-increasing}.
\end{problem}

\begin{problem}
Use Theorem~\ref{thm:polys} to study the unimodality, log-concavity and zeros of the polynomials of
$\G(\Z_n)$ for all $n$.
\end{problem}

\begin{problem}
Determine $\gamma(\G(L))$ for finite lattices whose maximal ideals are not all prime.
\end{problem}

\end{document}